\documentclass[a4paper]{amsart}
\usepackage[OT2,T1]{fontenc}
\usepackage{epic}
\theoremstyle{plain}
  \newtheorem{thm}{Theorem}[section]

  \newtheorem{conj}[thm]{Conjecture}
  
  \newtheorem{obs}[thm]{Observation}
\theoremstyle{definition}
  \newtheorem{defn}[thm]{Definition}
  
  \newtheorem{prob}[thm]{Problem}
\theoremstyle{remark}
  \newtheorem{rem}[thm]{Remark}
  \newtheorem*{ack}{Acknowledgments}
\DeclareMathAlphabet{\cyr}{OT2}{wncyr}{m}{n}
\newcommand{\Lob}{\operatorname{\cyr{L}}}

\newcommand{\Z}{\mathbb{Z}}
\newcommand{\C}{\mathbb{C}}

\newcommand{\Li}{\operatorname{Li}}
\newcommand{\cs}{\operatorname{cs}}
\newcommand{\CS}{\operatorname{CS}}
\newcommand{\Vol}{\operatorname{Vol}}

\DeclareMathOperator*{\olim}{o-lim}
\numberwithin{equation}{section}
\date{\today}
\allowdisplaybreaks
\begin{document}
\title[Calculations about the WRT invariants of closed three-manifolds]
{Optimistic calculations about \\
the Witten--Reshetikhin--Turaev invariants of
\\
closed three-manifolds obtained from the figure-eight knot
\\
by integral Dehn surgeries}
\author{Hitoshi Murakami}
\address{
Department of Mathematics,
Tokyo Institute of Technology,
Oh-okayama, Meguro, Tokyo 152-8551, Japan
}
\email{starshea@tky3.3web.ne.jp}
\begin{abstract}
I calculate {\em optimistically} asymptotic behaviors of the WRT $SU(2)$
invariants for the three-manifolds obtained from the figure-eight knot by
$p$-surgeries with $p=0,1,2,\dots,10$,
from which one can extract volumes and the Chern--Simons invariants of these
closed manifolds.
I conjecture that this also holds for general closed three-manifolds.
\end{abstract}
\maketitle
\section{Introduction}
In \cite{Kashaev:LETMP97} R.~Kashaev defined a link invariant by using quantum
dilogarithm and confirmed that his invariants grow exponentially with the growth
rates the hyperbolic volumes (times a constant) for three hyperbolic knots with
small numbers of crossings.
He also conjectured that this holds for every hyperbolic knot.
\par
J.~Murakami and I proved that Kashaev's link invariant is essentially (up to
normalization) the same as the Jones polynomial colored with $N$-dimensional
representation evaluated at the $N$-th root of unity.
Moreover we generalized Kashaev's conjecture to the following conjecture.
\begin{conj}[Volume Conjecture, \cite{Murakami/Murakami:volume}]
  Let $J_N(K)$ be the $N$-colored Jones polynomial of a knot $K$ evaluated
  at $\exp\left(\dfrac{2\pi\sqrt{-1}}{N}\right)$.
  Then
  \begin{equation*}
    \lim_{N\to\infty}\frac{\log|J_N(K)|}{N}=\frac{v_3}{2\pi}\Vert{K}\Vert,
  \end{equation*}
  where $\Vert{K}\Vert$ is the Gromov norm (or the simplicial volume) of the complement of $K$ and
  $v_3$ is the hyperbolic volume of the regular ideal tetrahedron.
\end{conj}
Recent developments toward the Volume Conjecture can be found in
\cite{Kashaev/Tirkkonen:1999,
Yokota:Murasugi70,
Yokota:volume,
Murakami:4_1,
Murakami/Murakami/Okamoto/Takata/Yokota:CS}.
\par
It is natural to ask whether a similar formula holds for closed three-manifolds
replacing the colored Jones polynomial with the Witten--Reshetikhin--Turaev
$SU(2)$ invariant associated with the $N$-th root of unity.
But an argument using Heegaard splitting and Topological Quantum Field
Theory tells us that the growth of the WRT invariant is a polynomial,
showing that a similar limit in the Volume Conjecture vanishes.
(After the first attempt of this work I learned this argument from D.~Thurston
and J.~Roberts; it was also pointed out by S.~Garoufalidis, V.~Turaev and
K.~Walker independently.)
\par
The aim of this article is to be optimistic to calculate (fake) limits of
the logarithms of the WRT invariants of closed three-manifolds obtained from
the figure-eight knot by Dehn surgeries with integral coefficients.
I will follow Kashaev's calculation in \cite{Kashaev:LETMP97} formally and
optimistically, and deduce an analytic function with integer parameter
corresponding to the surgery coefficient.
The function turns out to describe not only the (simplicial) volume but
also the Chern--Simons invariant of the manifold.
(J.~Murakami told me to look at the Chern--Simons invariants after my earlier
calculations.
See \cite{Murakami/Murakami/Okamoto/Takata/Yokota:CS} for a similar relation
between the Chern--Simons invariants and the colored Jones polynomials for knots
and links.)
\par
I do not know what these optimistic calculations mean.
But this is not a coincidence and there should be something behind it!
\begin{ack}
This article is for the proceedings of the workshop `Recent Progress toward the
Volume Conjecture' held at the International Institute for Advanced Study
from 14th to 17th March, 2000
(http://www.iias.or.jp/research/suuken/\linebreak[1]20000314eng.html),
supported by the Research Institute of Mathematical Sciences, Kyoto University.
I thank the IIAS for its hospitality and the RIMS for its financial support.

I am grateful to K.~Saito, who encouraged me to organize the workshop.
\par
Thanks are also due to the participants of the workshop, and to K.~Ichihara,
K.~Motegi, J.~Murakami, J.~Roberts,  M.~Teragaito, D.~Thurston and Y.~Yokota for
their helpful comments.
\end{ack}
\section{The Witten--Reshetikhin--Turaev invariant}
Let $J_n(K;t)$ be the colored Jones polynomial of a knot $K$ associated
with the $n$-dimensional representation of the Lie algebra $sl_2(\C)$.
We normalize $J_n(K;t)$ so that $J_2(K;t)$ is the Jones polynomial and
$J_n(O;t)=\dfrac{t^{n/2}-t^{-n/2}}{t^{1/2}-t^{-1/2}}$ with $O$ the unknot.
Due to T.~Le the colored Jones polynomial of the figure-eight knot $4_1$ is
\begin{equation*}
  J_n(4_1;t)=
  \sum_{m=0}^{n-1}
  \prod_{l=1}^{m}
  \left(t^{(n+l)/2}-t^{-(n+l)/2}\right)
  \left(t^{(n-l)/2}-t^{-(n-l)/2}\right).
\end{equation*}
(K.~Habiro obtained the same formula by a different technique.)
Let $M_p$ be the closed three-manifold obtained from $S^3$ by Dehn surgery
along the figure-eight knot with coefficient $p\in\Z$.
We denote by $\tau_N(M_p)$ the Witten--Reshetikhin--Turaev invariant of $M_p$
associated with the Lie group $sl_2(\C)$
and with level $N-2$.
Then from \cite[(1.9)]{Kirby/Melvin:INVEM91},
\begin{equation}
  \tau_{N}(M_p)
  =
  \sqrt{\dfrac{2}{N}}\sin\dfrac{\pi}{N}
  \exp\left(\dfrac{-3\pi\sqrt{-1}}{4}\right)q^{(3-p)/4}
  \sum_{n=1}^{N-1}[n]^2q^{pn^2/4}J_{n}(4_1;q)
\end{equation}
if $p>0$.
Here $q=\exp\left(\dfrac{2\pi\sqrt{-1}}{N}\right)$ and
$[n]=\dfrac{q^{n/2}-q^{-n/2}}{q^{1/2}-q^{-1/2}}$.
Note that our $J_n(K;q)$ is $[n]J_{K_0,n}$ with the notation in
\cite{Kirby/Melvin:INVEM91}.
Since
\begin{align*}
  [n]&
  \prod_{l=1}^{m}
  \left(q^{(n+l)/2}-q^{-(n+l)/2}\right)
  \left(q^{(n-l)/2}-q^{-(n-l)/2}\right)
  \\
  &=
  \dfrac{\prod_{l=-m}^{m}(q^{(n+l)/2}-q^{-(n+l)/2})}{q^{1/2}-q^{-1/2}}
  \\
  &=
  \dfrac{(1-q^{n-m})(1-q^{n-m+1})\cdots(1-q^{n+m})\times(-1)q^{-n(2m+1)/2}}
        {q^{1/2}-q^{-1/2}}
  \\
  &=
  \dfrac{-q^{-n(2m+1)/2}}{q^{1/2}-q^{-1/2}}\times\dfrac{(q)_{n+m}}{(q)_{n-m-1}}
\end{align*}
we have
\begin{align*}
  \tau_{N}(M_p)=
  P(N)
  \sum_{n=1}^{N-1}\sum_{m=0}^{n-1}
  \frac{(q)_n(q)_{n+m}}{(q)_{n-1}(q)_{n-m-1}}q^{n(pn/4-m)-n},
\end{align*}
where $(q)_k=(1-q)(1-q^2)\cdots(1-q^k)$ and $P(N)$ is a function of $N$ with
polynomial growth.
Note that this also holds for a non-positive $p$.
\section{An optimistic limit}
Suppose that we are given a function $S(N)$ of $N$ by the following
summation.
\begin{equation}\label{eq:sum}
  S(N)
  =
  P(N)
  \sum_{n_1,n_2,\dots,n_k}
  q^{Q+L}
  \prod_{a=1}^{\alpha}(q)^{\varepsilon_a}_{l_a},
\end{equation}
where $P(N)$ is a function of $N$ with polynomial growth,
$\varepsilon_a=\pm1$, $l_a$ and $L$ are linear functions of
$n_1,n_2,\dots,n_k$ (they do not depend on $N$ but may have constant terms),
$Q=\displaystyle\sum_{1\le i\le j\le k}r_{ij}n_in_j$,
and the summation runs over some range in
$\{(n_1,n_2,\dots,n_k) \mid 0\le n_i \le N-1\,(i=1,2,\dots,k)\}$.
Note that $q=\exp\left(\dfrac{2\pi\sqrt{-1}}{N}\right)$ and we regard $S$ as
a function of $N$ rather than $q$.
\par
Then {\em an optimistic limit} of $\dfrac{2\pi\sqrt{-1}\log{S(N)}}{N}$,
denoted by
$\displaystyle\olim_{N\to\infty}\dfrac{2\pi\sqrt{-1}\log{S(N)}}{N}$
is defined as follows.
\par
First we replace $S(N)/P(N)$ with the following iterated integral $I(N)$ along
some contours.
\begin{equation}\label{eq:integral}
  I(N):=
  \idotsint
  \exp
  \left[
    \dfrac{N}{2\pi\sqrt{-1}}
    V(z_1,z_2,\dots,z_k)
  \right]
  dz_1\,dz_2 \cdots dz_k.
\end{equation}
Here $V(z_1,z_2,\dots,z_k)$ is defined as follows.
Put
\begin{equation}\label{eq:V}
  \tilde{V}(z_1,z_2,\dots,z_k)
  :=
  -\sum_{a=1}^{\alpha}\varepsilon_a
     \left\{\Li_2\left(x_a\right)-\dfrac{\pi^2}{6}\right\}
  +\sum_{1\le i\le j\le k}r_{ij}\log z_i \log z_j,
\end{equation}
where $z_i=q^{n_i}$ and $x_a=q^{{l_a}'}$ with ${l_a}'$ the degree one term in
$l_a$, and $\Li_2(z)$ is Euler's dilogarithm defined by
\begin{equation*}
  \Li_2(z):=-\int_{0}^{z}\dfrac{\log(1-u)}{u}\,du.
\end{equation*}
\par
Next we consider the following system of partial differential equations:
\begin{equation}\label{eq:differential}
  \dfrac{\partial\,\tilde{V}(z_1,z_2,\dots,z_k)}{\partial\,z_i}=0
  \quad(i=1,2,\dots,k).
\end{equation}
Since
\begin{equation*}
  \dfrac{\partial\,\tilde{V}(z_1,z_2,\dots,z_k)}{\partial\,z_i}
  =
  \sum_{a=1}^{\alpha}\varepsilon_al_{ai}\dfrac{\log(1-x_a)}{z_i}
  +
  \sum_{j=1}^{k}r_{ij}\dfrac{\log z_j}{z_i}+r_{ii}\dfrac{\log z_i}{z_i}
\end{equation*}
with ${l_a}'=\sum_{i=1}^{k}l_{ai}z_i$,
\eqref{eq:differential} implies the following algebraic equations.
\begin{equation}\label{eq:algebraic}
  {z_i}^{r_{ii}}\prod_{j=1}^{k}{z_j}^{r_{ij}}
  \prod_{a=1}^{\alpha}(1-x_a)^{\varepsilon_a l_{ai}}
  =1
  \quad(i=1,2,\dots,k).
\end{equation}
Let $(\zeta_1,\zeta_2,\dots,\zeta_k)$ be a solution to \eqref{eq:algebraic}.
\begin{defn}[optimistic limit]
  We put
  \begin{equation}\label{eq:definition}
  V(\zeta_1,\zeta_2,\dots,\zeta_k):=
  \tilde{V}(\zeta_1,\zeta_2,\dots,\zeta_k)
  +2\pi\sqrt{-1}\left(\sum_{j=1}^{k}c_j\log\zeta_j\right)
  \end{equation}
  and call it {\em an optimistic limit} of $\dfrac{2\pi\sqrt{-1}\log{S(N)}}{N}$
  as $N$ goes to the infinity.
  It is denoted by
  $\displaystyle\olim_{N\to\infty}\dfrac{2\pi\sqrt{-1}\log{S(N)}}{N}$.
  Here $c_i$ is chosen so that
  $\dfrac{\partial\tilde{V}(\zeta_1,\zeta_2,\dots,\zeta_k)}{\partial z_i}
  +2\pi\sqrt{-1}\sum_{j=1}^{k}\dfrac{c_j}{\zeta_j}=0$
  for every $i$.
\end{defn}
\begin{rem}
The term $-\dfrac{\pi^2}{6}$ in \eqref{eq:V} appears so that
$V(1,1,\dots,1)=0$ since $\Li_2(1)=\dfrac{\pi^2}{6}$ (see for example
\cite[(1.5)]{Kirillov:dilog}).
(I learned this from \cite[\S5]{Jun:MSJ2000}.)
\end{rem}
\begin{rem}
Note that $V$ and $\tilde{V}$ satisfy the same algebraic equations
\eqref{eq:algebraic} but different partial differential equations
\eqref{eq:differential}
and that the extra terms in \eqref{eq:definition} are
necessary to choose an appropriate branch since $\Li_2$ and $\log$ are
multivalued functions.
(I learned this from T.~Takata.)
\end{rem}
\begin{rem}
An optimistic limit is {\em not} well defined.
There are many ambiguities both in choosing $I(N)$ (I did not say anything
about the range of the summation in $S(N)$ and the contours in $I(N)$)
and in choosing $(\zeta_1,\zeta_2,\dots,\zeta_k)$.
\end{rem}
\begin{rem}
Following Kashaev \cite{Kashaev:LETMP97}, the behavior of $S(N)$ for
large $N$ may be approximated by $P(N)I(N)$ with suitably chosen contours.
Moreover by using the saddle point method (see for example
\cite[\S7.2]{Marsden/Hoffman:Complex_Analysis})
we see that $I(N)$ (and $S(N)$) behaves like
$\displaystyle\exp\left(\dfrac{N}{2\pi\sqrt{-1}}\olim_{N\to\infty}
 \dfrac{2\pi\sqrt{-1}\log{S(N)}}{N}\right)$
for large $N$ if we choose a solution $(\zeta_1,\zeta_2,\dots,\zeta_k)$
suitably.
Even if this is not true I expect that there is a relation between
an optimistic limit and the asymptotic behavior of $S(N)$.
\end{rem}
\par
\section{Dehn surgery along the figure-eight knot}
Put $k:=2$, $n_1:=n$, $n_2:=m$, $z:=z_1$, $w:=z_2$, $Q:=pn^2/4-mn$, $L:=-n/2$,
$\alpha:=4$, $l_1:=n$, $\varepsilon_1:=1$, $l_2:=n-1$, $\varepsilon_2:=-1$,
$l_3:=n+m$, $\varepsilon_3:=1$, $l_4:=n-m-1$, and $\varepsilon_4:=-1$
in \eqref{eq:sum}.
Note that
$x_1=z$, $x_2=z$, $x_3=zw$, $x_4=zw^{-1}$,
$r_{11}=p/4$, $r_{12}=-1$, $r_{22}=0$,
$l_{11}=1$, $l_{12}=0$, $l_{21}=1$, $l_{22}=0$,
$l_{31}=1$, $l_{32}=1$, $l_{41}=1$, $l_{42}=-1$.
Then
\begin{equation*}
  \tilde{V}(z,w):=
  -\Li_2(zw)+\Li_2\left(\dfrac{z}{w}\right)+
  \dfrac{p}{4}(\log z)^2- \log z\log w
\end{equation*}
and
\begin{equation*}
\begin{cases}
  \dfrac{\partial\,\tilde{V}}{\partial\,z}
  &=
  \dfrac{1}{z}
  \left\{
    \log z^{p/2}+\log\left(\dfrac{1-zw}{w-z}\right)
  \right\},
  \\[5mm]
  \dfrac{\partial\,\tilde{V}}{\partial\,w}
  &=
  \dfrac{1}{w}
  \log\dfrac{(1-zw)(w-z)}{zw}.
\end{cases}
\end{equation*}
Therefore \eqref{eq:algebraic} turns out to be
\begin{align}
  &
  \begin{cases}
    z^{p/2}(1-zw)=w-z,
    \\[5mm]
    (1-zw)(w-z)=zw,
  \end{cases}
  \label{eq:algebraic_fig8_original}
  \\ \notag
  \intertext{from which we have}
  &
  \begin{cases}
  w=\dfrac{z+z^{p/2}}{z^{p/2}z+1},
  \\[5mm]
  z^2-\left(\dfrac{z+z^{p/2}}{z^{p/2}z+1}+1
  +\dfrac{z^{p/2}z+1}{z+z^{p/2}}\right)z+1=0.
  \end{cases}
  \label{eq:algebraic_fig8}
\end{align}
\par
\begin{rem}
Since the second equation of \eqref{eq:algebraic_fig8} is symmetric
with respect to $z$ and $z^{-1}$, if $\zeta$ is a solution to it then
so is $\zeta^{-1}$.
(This may be caused by the amphicheirality of the figure-eight knot.)
Clearly $\overline{\zeta}$, the complex conjugate of $\zeta$, also satisfies it.
Therefore if $(\zeta,\omega)$ is a solution to \eqref{eq:algebraic_fig8}
then so are $(\overline{\zeta},\overline{\omega})$, $(\zeta^{-1},\omega)$,
and $(\overline{\zeta}^{-1},\overline{\omega})$.
\end{rem}
\par
I will show calculations for $p=0,1,\dots,10$.
\subsection{$6$-surgery along the figure-eight knot}
I will describe the case where $p=6$ in detail.
Note that $M_6$ is hyperbolic.
In this case there are the following six solutions to
\eqref{eq:algebraic_fig8} due to MAPLE V:
\begin{equation*}
  (\zeta_1,\omega_1),\,
  (\zeta_2,\omega_2),\,
  ({\zeta_1}^{-1},\omega_1),\,
  (\overline{\zeta_2},\overline{\omega_2}),\,
  ({\zeta_2}^{-1},\omega_2),\,
  ({\overline{\zeta_2}}^{-1},\overline{\omega_2}),
\end{equation*}
where
\begin{align*}
  (\zeta_1,\omega_1)
  &=
  (-0.8294835410-0.5585311587\sqrt{-1},
  -2.205569430\phantom{2}-0.3703811357\times10^{-9}\sqrt{-1})
  \\ \intertext{and}
  (\zeta_2,\omega_2)
  &=
  (\phantom{-}0.3679390314-0.4972675889\sqrt{-1},
  \phantom{-}0.1027847152-0.6654569513\sqrt{-1}).
\end{align*}
Note that $\left|\zeta_1\right|=1$ and so
$\overline{\zeta_1}={\zeta_1}^{-1}$ and ${\overline{\zeta_1}}^{-1}=\zeta_1$.
\par
The partial derivatives
$\dfrac{\partial\,\tilde{V}(\zeta,\omega)}{\partial\,z}$
and $\dfrac{\partial\,\tilde{V}(\zeta,\omega)}{\partial\,w}$
for $(\zeta_1,\omega_1)$ and $(\zeta_2,\omega_2)$ are as follows.
\begin{align*}
  &\begin{cases}
    \dfrac{\partial\,\tilde{V}}{\partial\,z}(\zeta_1,\omega_1)
    &=
    0.424142903\times10^{-10} -6.283185309\sqrt{-1},
    \\[5mm]
    \dfrac{\partial\,\tilde{V}}{\partial\,w}(\zeta_1,\omega_1)
    &=
    0.2205569430\times10^{-9} -0.2205569430\times10^{-9}\sqrt{-1},
  \end{cases}
  \\ \intertext{and}
  &\begin{cases}
    \dfrac{\partial\,\tilde{V}}{\partial\,z}(\zeta_2,\omega_2)
    &=
    0.3868858795\times10^{-9} +0.5171193081\times10^{-9}\sqrt{-1},
    \\[5mm]
    \dfrac{\partial\,\tilde{V}}{\partial\,w}(\zeta_2,\omega_2)
    &=
    0.3623902007\times10^{-10} -0.6757354228\times10^{-9}\sqrt{-1}.
  \end{cases}
\end{align*}
Therefore we have
\begin{align*}
  V(\zeta_1,\omega_1)
  &=
  \tilde{V}(\zeta_1,\omega_1)+2\pi\sqrt{-1}\log\zeta_1,
  \\
  V(\zeta_2,\omega_2)
  &=
  \tilde{V}(\zeta_2,\omega_2),
  \\
  V({\zeta_1}^{-1},\omega_1)
  &=
  V(\zeta_1,\omega_1),
  \\
  V(\overline{\zeta_2},\overline{\omega_2})
  &=
  \overline{V(\zeta_2,\omega_2)},
  \\
  V({\zeta_2}^{-1},\omega_2)
  &=
  V(\zeta_2,\omega_2),
  \\
  V({\overline{\zeta_2}}^{-1},\overline{\omega_2})
  &=
  \overline{V(\zeta_2,\omega_2)},
\end{align*}
with
\begin{align*}
  V(\zeta_1,\omega_1)&=13.76750570 +0.1\times10^{-8}\sqrt{-1}
  \\ \intertext{and}
  V(\zeta_2,\omega_2)&=1.340917487 +1.284485301\sqrt{-1}
\end{align*}
by MAPLE V.
\par
So there are three optimistic limits (up to 10 digits);
$V_1:=13.76750570$, $V_2:=1.340917487+1.284485301\sqrt{-1}$,
and $\overline{V_2}$.
By using SnapPea \cite{Weeks:SnapPea}, we calculate
$\Vol(M_6)=1.2844853$ and $\cs(M_6)=0.0679316734799$
and so we can write
\begin{equation*}
  V_2=\CS(M_6)+\Vol(M_6)\sqrt{-1}.
\end{equation*}
since $0.0679316734799\times2\pi^2=1.34091748750\dots$.
Here $\Vol(M)$ and $\cs(M)$ are the volume and the Chern--Simons invariant
\cite{Chern/Simons:ANNMA274} of a closed hyperbolic three-manifold $M$
respectively, and $\CS(M):=2\pi^2\cs(M)$.
\subsection{$5,7,8,9,10$-surgeries along the figure-eight knot}
Similar calculations using MAPLE V for $p=5,7,8,9,10$ give the following.
Note that $M_p$ is also hyperbolic in this case.
\begin{obs}\label{obs}
  For $p=5,6,7,8,9,10$ there is an optimistic limit $V(\zeta_p,\omega_p)$
  of $\dfrac{2\pi\sqrt{-1}\log\tau_N(M_p)}{N}$ such that
  \begin{equation*}
    V(\zeta_p,\omega_p)=\CS(M_p)+\Vol(M_p)\sqrt{-1}
  \end{equation*}
  up to several digits.
  Here
  \begin{equation*}
    V(z,w)=-\Li_2(zw)+\Li_2(z/w)+\dfrac{p}{4}\left(\log{z}\right)^2
           -\log{z}\log{w}
  \end{equation*}
  and
  \begin{align*}
    (\zeta_5,\omega_5)&=(0.1979823656-0.4438341209\sqrt{-1},\,
                         0.007552359501-0.5131157955\sqrt{-1}),
    \\
    (\zeta_6,\omega_6)&=(0.3679390314-0.4972675889\sqrt{-1},\,
                         0.1027847152\phantom{00}-0.6654569513\sqrt{-1}),
    \\
    (\zeta_7,\omega_7)&=(0.4855046904-0.5042960525\sqrt{-1},\,
                         0.1761405059\phantom{00}-0.7455559248\sqrt{-1}),
    \\
    (\zeta_8,\omega_8)&=(0.5730134132-0.4940983127\sqrt{-1},\,
                         0.2327856161\phantom{00}-0.7925519927\sqrt{-1}),
    \\
    (\zeta_9,\omega_9)&=(0.6404276706-0.4765868179\sqrt{-1},\,
                         0.2769632324\phantom{00}-0.8216401587\sqrt{-1}),
    \\
    (\zeta_{10},\omega_{10})&=(0.6935298015-0.4561607978\sqrt{-1},\,
                               0.3118108269\phantom{00}-0.8402398912\sqrt{-1}).
  \end{align*}
\end{obs}
\begin{rem}
  For $p=4n+2$ with $n=1,2,\dots,100$, we can observe the same result.
  See Figures~1 and 2 for $(\zeta_p,\omega_p)$.
\end{rem}
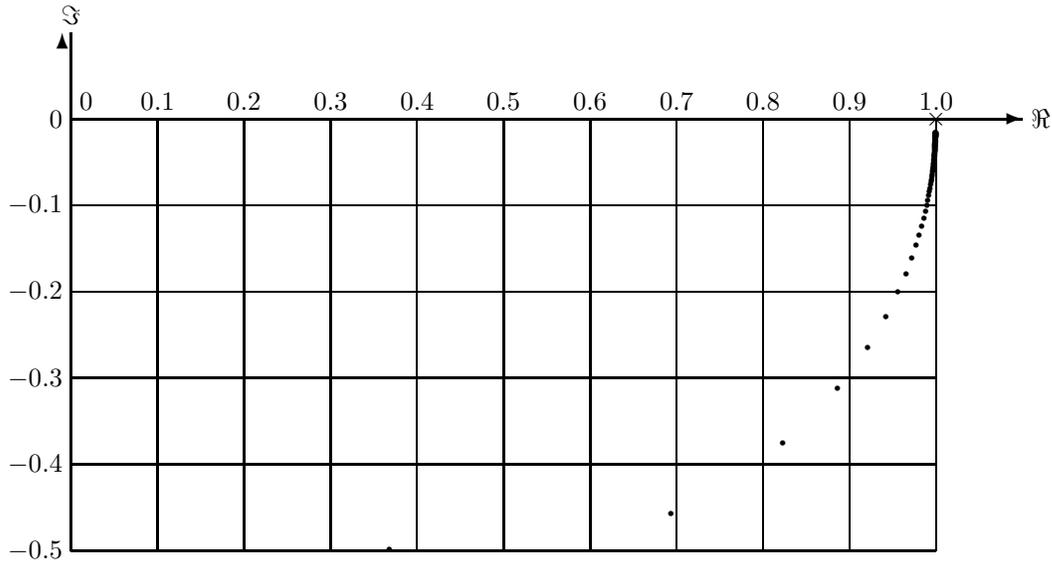
\begin{figure}[h]
\setlength{\unitlength}{11.5cm}
\begin{picture}(1.1,0.6)(0,-0.5)
\thicklines
\put(0,0.0){\vector(1,0){1.1}}
\put(0,-0.5){\line(0,1){0.6}}
\put(-0.01,0.1){\vector(0,1){0}}
\thinlines
\put(0,-0.5){\grid(1,0.5)(0.1,0.1)}
\put(1,-0.5){\line(0,1){0.5}}
\put(0,0){\line(1,0){1}}
\put(0.01,0.01){\makebox(0,0)[bl]{$0$}}
\put(0.1,0.01){\makebox(0,0)[b]{$0.1$}}
\put(0.2,0.01){\makebox(0,0)[b]{$0.2$}}
\put(0.3,0.01){\makebox(0,0)[b]{$0.3$}}
\put(0.4,0.01){\makebox(0,0)[b]{$0.4$}}
\put(0.5,0.01){\makebox(0,0)[b]{$0.5$}}
\put(0.6,0.01){\makebox(0,0)[b]{$0.6$}}
\put(0.7,0.01){\makebox(0,0)[b]{$0.7$}}
\put(0.8,0.01){\makebox(0,0)[b]{$0.8$}}
\put(0.9,0.01){\makebox(0,0)[b]{$0.9$}}
\put(1.0,0.01){\makebox(0,0)[b]{$1.0$}}
\put(-0.01,-0.0){\makebox(0,0)[r]{$0$}}
\put(-0.01,-0.1){\makebox(0,0)[r]{$-0.1$}}
\put(-0.01,-0.2){\makebox(0,0)[r]{$-0.2$}}
\put(-0.01,-0.3){\makebox(0,0)[r]{$-0.3$}}
\put(-0.01,-0.4){\makebox(0,0)[r]{$-0.4$}}
\put(-0.01,-0.5){\makebox(0,0)[r]{$-0.5$}}
\put(1.11,0){\makebox(0,0)[l]{$\Re$}}
\put(0,0.11){\makebox(0,0)[b]{$\Im$}}
\put(1,0){\makebox(0,0){$\times$}}
\put(.3679390314, -.4972675889){\circle*{0.007}}
\put(.6935298015, -.4561607978){\circle*{0.007}}
\put(.8226587769, -.3752413832){\circle*{0.007}}
\put(.8857985824, -.3115833941){\circle*{0.007}}
\put(.9208304197, -.2641278732){\circle*{0.007}}
\put(.9420942101, -.2282907736){\circle*{0.007}}
\put(.9558977577, -.2005785405){\circle*{0.007}}
\put(.9653358234, -.1786371202){\circle*{0.007}}
\put(.9720602272, -.1608939285){\circle*{0.007}}
\put(.9770137446, -.1462800976){\circle*{0.007}}
\put(.9807646131, -.1340518161){\circle*{0.007}}
\put(.9836711388, -.1236789210){\circle*{0.007}}
\put(.9859679887, -.1147749309){\circle*{0.007}}
\put(.9878138992, -.1070522647){\circle*{0.007}}
\put(.9893192712, -.1002929384){\circle*{0.007}}
\put(.9905627676, -.09432895235){\circle*{0.007}}
\put(.9916016380, -.08902888119){\circle*{0.007}}
\put(.9924783357, -.08428851949){\circle*{0.007}}
\put(.9932248715, -.08002423431){\circle*{0.007}}
\put(.9938657450, -.07616816089){\circle*{0.007}}
\put(.9944199614, -.07266467592){\circle*{0.007}}
\put(.9949024436, -.06946777130){\circle*{0.007}}
\put(.9953250401, -.06653907271){\circle*{0.007}}
\put(.9956972536, -.06384632642){\circle*{0.007}}
\put(.9960267761, -.06136223079){\circle*{0.007}}
\put(.9963198878, -.05906352455){\circle*{0.007}}
\put(.9965817561, -.05693026868){\circle*{0.007}}
\put(.9968166643, -.05494527556){\circle*{0.007}}
\put(.9970281866, -.05309365171){\circle*{0.007}}
\put(.9972193246, -.05136242840){\circle*{0.007}}
\put(.9973926136, -.04974026128){\circle*{0.007}}
\put(.9975502070, -.04821718436){\circle*{0.007}}
\put(.9976939433, -.04678440738){\circle*{0.007}}
\put(.9978253996, -.04543414771){\circle*{0.007}}
\put(.9979459350, -.04415949032){\circle*{0.007}}
\put(.9980567258, -.04295427038){\circle*{0.007}}
\put(.9981587940, -.04181297437){\circle*{0.007}}
\put(.9982530315, -.04073065647){\circle*{0.007}}
\put(.9983402188, -.03970286743){\circle*{0.007}}
\put(.9984210416, -.03872559386){\circle*{0.007}}
\put(.9984961043, -.03779520622){\circle*{0.007}}
\put(.9985659414, -.03690841407){\circle*{0.007}}
\put(.9986310267, -.03606222742){\circle*{0.007}}
\put(.9986917814, -.03525392325){\circle*{0.007}}
\put(.9987485813, -.03448101639){\circle*{0.007}}
\put(.9988017622, -.03374123414){\circle*{0.007}}
\put(.9988516250, -.03303249410){\circle*{0.007}}
\put(.9988984400, -.03235288460){\circle*{0.007}}
\put(.9989424505, -.03170064766){\circle*{0.007}}
\put(.9989838762, -.03107416382){\circle*{0.007}}
\put(.9990229154, -.03047193876){\circle*{0.007}}
\put(.9990597479, -.02989259149){\circle*{0.007}}
\put(.9990945370, -.02933484380){\circle*{0.007}}
\put(.9991274309, -.02879751087){\circle*{0.007}}
\put(.9991585650, -.02827949291){\circle*{0.007}}
\put(.9991880623, -.02777976768){\circle*{0.007}}
\put(.9992160357, -.02729738376){\circle*{0.007}}
\put(.9992425884, -.02683145454){\circle*{0.007}}
\put(.9992678149, -.02638115281){\circle*{0.007}}
\put(.9992918020, -.02594570582){\circle*{0.007}}
\put(.9993146297, -.02552439095){\circle*{0.007}}
\put(.9993363714, -.02511653164){\circle*{0.007}}
\put(.9993570950, -.02472149383){\circle*{0.007}}
\put(.9993768629, -.02433868264){\circle*{0.007}}
\put(.9993957332, -.02396753943){\circle*{0.007}}
\put(.9994137592, -.02360753905){\circle*{0.007}}
\put(.9994309907, -.02325818739){\circle*{0.007}}
\put(.9994474736, -.02291901910){\circle*{0.007}}
\put(.9994632506, -.02258959554){\circle*{0.007}}
\put(.9994783615, -.02226950292){\circle*{0.007}}
\put(.9994928433, -.02195835051){\circle*{0.007}}
\put(.9995067304, -.02165576909){\circle*{0.007}}
\put(.9995200549, -.02136140951){\circle*{0.007}}
\put(.9995328467, -.02107494128){\circle*{0.007}}
\put(.9995451340, -.02079605141){\circle*{0.007}}
\put(.9995569429, -.02052444324){\circle*{0.007}}
\put(.9995682979, -.02025983537){\circle*{0.007}}
\put(.9995792219, -.02000196073){\circle*{0.007}}
\put(.9995897366, -.01975056563){\circle*{0.007}}
\put(.9995998620, -.01950540900){\circle*{0.007}}
\put(.9996096172, -.01926626155){\circle*{0.007}}
\put(.9996190201, -.01903290510){\circle*{0.007}}
\put(.9996280873, -.01880513188){\circle*{0.007}}
\put(.9996368346, -.01858274397){\circle*{0.007}}
\put(.9996452770, -.01836555267){\circle*{0.007}}
\put(.9996534284, -.01815337800){\circle*{0.007}}
\put(.9996613020, -.01794604820){\circle*{0.007}}
\put(.9996689104, -.01774339923){\circle*{0.007}}
\put(.9996762653, -.01754527442){\circle*{0.007}}
\put(.9996833779, -.01735152396){\circle*{0.007}}
\put(.9996902586, -.01716200462){\circle*{0.007}}
\put(.9996969174, -.01697657933){\circle*{0.007}}
\put(.9997033639, -.01679511689){\circle*{0.007}}
\put(.9997096068, -.01661749164){\circle*{0.007}}
\put(.9997156548, -.01644358318){\circle*{0.007}}
\put(.9997215157, -.01627327607){\circle*{0.007}}
\put(.9997271974, -.01610645964){\circle*{0.007}}
\put(.9997327069, -.01594302768){\circle*{0.007}}
\put(.9997380512, -.01578287824){\circle*{0.007}}
\put(.9997432368, -.01562591345){\circle*{0.007}}
\end{picture}
\caption{$\zeta_p$ is indicated by a dot for $p=4n+2$ ($n=1,2,\dots,100$).
         The dots approach $\zeta_{\infty}=1$ indicated by $\times$ from
         the left when $n$ increases.}
\end{figure}
\begin{figure}[h]
\setlength{\unitlength}{11.5cm}
\begin{picture}(1.1,0.6)(0,-1)
\thicklines
\put(0,-0.5){\vector(1,0){1.1}}
\put(0,-1){\line(0,1){0.6}}
\put(-0.01,-0.4){\vector(0,1){0}}
\thinlines
\put(0,-1){\grid(1,0.5)(0.1,0.1)}
\put(1,-1){\line(0,1){0.5}}
\put(0,-0.5){\line(1,0){1}}
\put(0.01,-0.49){\makebox(0,0)[bl]{$0$}}
\put(0.1,-0.49){\makebox(0,0)[b]{$0.1$}}
\put(0.2,-0.49){\makebox(0,0)[b]{$0.2$}}
\put(0.3,-0.49){\makebox(0,0)[b]{$0.3$}}
\put(0.4,-0.49){\makebox(0,0)[b]{$0.4$}}
\put(0.5,-0.49){\makebox(0,0)[b]{$0.5$}}
\put(0.6,-0.49){\makebox(0,0)[b]{$0.6$}}
\put(0.7,-0.49){\makebox(0,0)[b]{$0.7$}}
\put(0.8,-0.49){\makebox(0,0)[b]{$0.8$}}
\put(0.9,-0.49){\makebox(0,0)[b]{$0.9$}}
\put(1.0,-0.49){\makebox(0,0)[b]{$1.0$}}
\put(-0.01,-0.5){\makebox(0,0)[r]{$-0.5$}}
\put(-0.01,-0.6){\makebox(0,0)[r]{$-0.6$}}
\put(-0.01,-0.7){\makebox(0,0)[r]{$-0.7$}}
\put(-0.01,-0.8){\makebox(0,0)[r]{$-0.8$}}
\put(-0.01,-0.9){\makebox(0,0)[r]{$-0.9$}}
\put(-0.01,-1.0){\makebox(0,0)[r]{$-1.0$}}
\put(1.11,-0.5){\makebox(0,0)[l]{$\Re$}}
\put(0,-0.39){\makebox(0,0)[b]{$\Im$}}
\put(0.5,-.8660254040){\makebox(0,0){$\times$}}
\put(.1027847152, -.6654569513){\circle*{0.007}}
\put(.3118108269, -.8402398912){\circle*{0.007}}
\put(.3953991688, -.8693741692){\circle*{0.007}}
\put(.4348230327, -.8748318308){\circle*{0.007}}
\put(.4559533439, -.8750705282){\circle*{0.007}}
\put(.4684153953, -.8741238229){\circle*{0.007}}
\put(.4763180163, -.8730162779){\circle*{0.007}}
\put(.4816190375, -.8720115502){\circle*{0.007}}
\put(.4853367701, -.8711609688){\circle*{0.007}}
\put(.4880395674, -.8704561515){\circle*{0.007}}
\put(.4900634817, -.8698745703){\circle*{0.007}}
\put(.4916169363, -.8693933227){\circle*{0.007}}
\put(.4928344808, -.8689927256){\circle*{0.007}}
\put(.4938059979, -.8686569400){\circle*{0.007}}
\put(.4945933202, -.8683733648){\circle*{0.007}}
\put(.4952400593, -.8681321238){\circle*{0.007}}
\put(.4957776956, -.8679254635){\circle*{0.007}}
\put(.4962294000, -.8677472309){\circle*{0.007}}
\put(.4966124890, -.8675925655){\circle*{0.007}}
\put(.4969401782, -.8674575467){\circle*{0.007}}
\put(.4972225998, -.8673390603){\circle*{0.007}}
\put(.4974677526, -.8672345147){\circle*{0.007}}
\put(.4976818797, -.8671418521){\circle*{0.007}}
\put(.4978699685, -.8670593756){\circle*{0.007}}
\put(.4980361255, -.8669856175){\circle*{0.007}}
\put(.4981835806, -.8669194427){\circle*{0.007}}
\put(.4983150666, -.8668598161){\circle*{0.007}}
\put(.4984327844, -.8668059400){\circle*{0.007}}
\put(.4985386140, -.8667570767){\circle*{0.007}}
\put(.4986340690, -.8667126512){\circle*{0.007}}
\put(.4987204856, -.8666721208){\circle*{0.007}}
\put(.4987989877, -.8666350477){\circle*{0.007}}
\put(.4988704711, -.8666010685){\circle*{0.007}}
\put(.4989357534, -.8665698510){\circle*{0.007}}
\put(.4989955329, -.8665410997){\circle*{0.007}}
\put(.4990503897, -.8665145697){\circle*{0.007}}
\put(.4991009677, -.8664899740){\circle*{0.007}}
\put(.4991475328, -.8664672409){\circle*{0.007}}
\put(.4991905562, -.8664461402){\circle*{0.007}}
\put(.4992304414, -.8664264864){\circle*{0.007}}
\put(.4992675096, -.8664081369){\circle*{0.007}}
\put(.4993019041, -.8663910559){\circle*{0.007}}
\put(.4993338843, -.8663751228){\circle*{0.007}}
\put(.4993637970, -.8663601453){\circle*{0.007}}
\put(.4993917127, -.8663461295){\circle*{0.007}}
\put(.4994178190, -.8663329975){\circle*{0.007}}
\put(.4994422805, -.8663206269){\circle*{0.007}}
\put(.4994652590, -.8663089905){\circle*{0.007}}
\put(.4994869114, -.8662979823){\circle*{0.007}}
\put(.4995071258, -.8662876881){\circle*{0.007}}
\put(.4995262399, -.8662779316){\circle*{0.007}}
\put(.4995443313, -.8662686692){\circle*{0.007}}
\put(.4995612936, -.8662599714){\circle*{0.007}}
\put(.4995774652, -.8662516407){\circle*{0.007}}
\put(.4995925183, -.8662439062){\circle*{0.007}}
\put(.4996070205, -.8662363957){\circle*{0.007}}
\put(.4996207193, -.8662292979){\circle*{0.007}}
\put(.4996336079, -.8662226350){\circle*{0.007}}
\put(.4996458644, -.8662162522){\circle*{0.007}}
\put(.4996576576, -.8662101157){\circle*{0.007}}
\put(.4996687183, -.8662043494){\circle*{0.007}}
\put(.4996793441, -.8661987962){\circle*{0.007}}
\put(.4996893304, -.8661935727){\circle*{0.007}}
\put(.4996991188, -.8661884380){\circle*{0.007}}
\put(.4997081689, -.8661837088){\circle*{0.007}}
\put(.4997170428, -.8661790343){\circle*{0.007}}
\put(.4997253067, -.8661746832){\circle*{0.007}}
\put(.4997332630, -.8661704891){\circle*{0.007}}
\put(.4997410350, -.8661664006){\circle*{0.007}}
\put(.4997484293, -.8661624757){\circle*{0.007}}
\put(.4997553982, -.8661587891){\circle*{0.007}}
\put(.4997620768, -.8661552785){\circle*{0.007}}
\put(.4997684933, -.8661518455){\circle*{0.007}}
\put(.4997749405, -.8661483994){\circle*{0.007}}
\put(.4997808543, -.8661452756){\circle*{0.007}}
\put(.4997864949, -.8661422704){\circle*{0.007}}
\put(.4997919338, -.8661393871){\circle*{0.007}}
\put(.4997973763, -.8661364489){\circle*{0.007}}
\put(.4998023090, -.8661338648){\circle*{0.007}}
\put(.4998073810, -.8661311377){\circle*{0.007}}
\put(.4998121819, -.8661285576){\circle*{0.007}}
\put(.4998164451, -.8661262711){\circle*{0.007}}
\put(.4998208238, -.8661239521){\circle*{0.007}}
\put(.4998252447, -.8661215616){\circle*{0.007}}
\put(.4998292234, -.8661194454){\circle*{0.007}}
\put(.4998331920, -.8661173188){\circle*{0.007}}
\put(.4998372171, -.8661150921){\circle*{0.007}}
\put(.4998408547, -.8661131825){\circle*{0.007}}
\put(.4998443673, -.8661112457){\circle*{0.007}}
\put(.4998475511, -.8661095758){\circle*{0.007}}
\put(.4998509433, -.8661077315){\circle*{0.007}}
\put(.4998544358, -.8661058429){\circle*{0.007}}
\put(.4998572556, -.8661043580){\circle*{0.007}}
\put(.4998605180, -.8661025692){\circle*{0.007}}
\put(.4998631373, -.8661011327){\circle*{0.007}}
\put(.4998663395, -.8660994108){\circle*{0.007}}
\put(.4998686977, -.8660981458){\circle*{0.007}}
\put(.4998714677, -.8660966058){\circle*{0.007}}
\put(.4998740807, -.8660952274){\circle*{0.007}}
\put(.4998767005, -.8660937909){\circle*{0.007}}
\end{picture}
\caption{$\omega_p$ is indicated by a dot for $p=4n+2$ ($n=1,2,\dots,100$).
         The dots approach
         $\omega_{\infty}=\exp\left(-\dfrac{\pi\sqrt{-1}}{3}\right)$
         indicated by $\times$ from the left when $n$ increases.}
\end{figure}
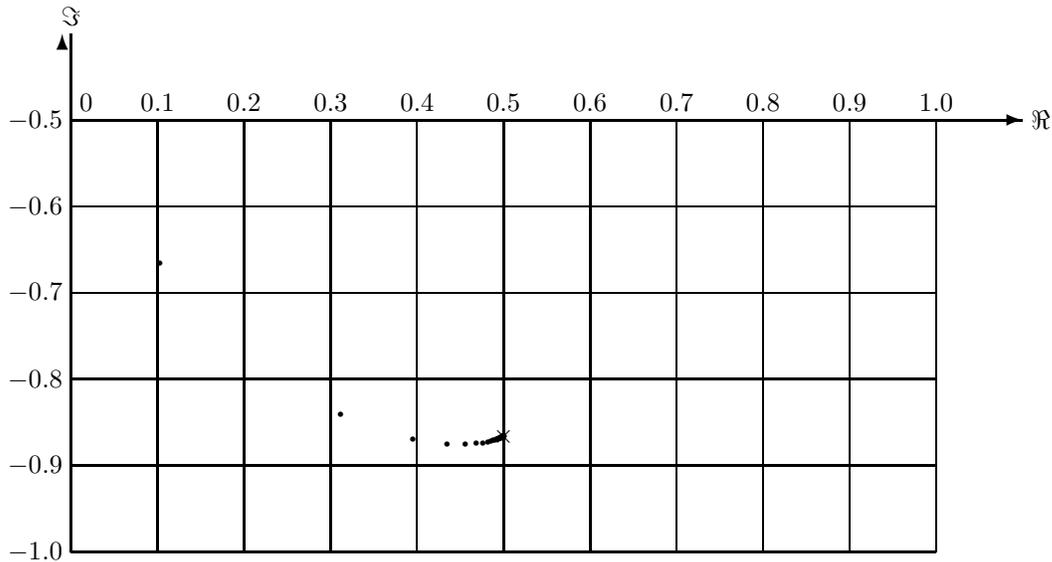
\subsection{$1,2,3$-surgeries along the figure-eight knot}
Next we consider the case where $p=1,2$ and $3$.
Note that the manifold $M_p$ for $p=1,2,3$ is a Seifert fibered space.
See \cite[p.~95]{Kirby:problems} for details, which was informed by
K.~Ichihara.
\par
In this case, the (simplicial) volume of $M_p$ is zero but SnapPea tells us
that it has the non-trivial Chern--Simons invariant.
MAPLE V shows that the same observation as Observation~\ref{obs} holds
with
\begin{align*}
  (\zeta_1,\omega_1)&=(0.3738178762,0.8019377355),
  \\
  (\zeta_2,\omega_2)&=(0.346014339\phantom{0},0.6180339884),
  \\
  (\zeta_3,\omega_3)&=(0.2819716801,0.4142135623).
\end{align*}
\subsection{$0$-surgery of the figure-eight knot}
In this case $M_0$ is a torus bundle over a circle.
For a detail see \cite[p.~95]{Kirby:problems} again.
Both the volume and the Chern--Simons invariant vanish in this case and
Observation~\ref{obs} holds putting $(\zeta_0,\omega_0)=(0.381966011,1)$.
\subsection{$4$-surgery of the figure-eight knot}
It is known that $M_4$ is toroidal and can be obtained by gluing the
twisted $I$-bundle over the Klein bottle and the complement of the trefoil,
which I learned from K.~Motegi and M.~Teragaito.
The gluing map can be found in \cite[p.~95]{Kirby:problems} again.
Therefore $\Vol(M_4)=0$.
(Here I use $\Vol$ for $v_3$ times the simplicial volume.)
Computation by using MAPLE V told us that
if we put $(\zeta_4,\omega_4)=(-1,-.381966011)$ then
$V(\zeta_4,\omega_4)=1.973920880+0.1\times10^{-8}\sqrt{-1}$.
Note that $1.973920880=2\pi^2\times0.09999999995$.
I guess this is the sum of the Chern--Simons invariants of
the two pieces (with suitably chosen metrics),
which one might regard as the Chern--Simons invariant of
the toroidal manifold $M_4$.
\subsection{$\infty$-surgery of the figure-eight knot}
Put $\zeta_{\infty}:=\exp\left(-\dfrac{2\pi\sqrt{-1}}{p}\right)$ and
$\omega_{\infty}:=\exp\left(-\dfrac{\pi\sqrt{-1}}{3}\right)$.
Then the left hand side minus the right hand side of the first equation
in \eqref{eq:algebraic_fig8_original} is
$(\zeta_{\infty}-1)(\omega_{\infty}+1)$ since ${\zeta_{\infty}}^{p/2}=-1$.
That of the second equation is $\omega_{\infty}(\zeta_{\infty}-1)^2$
since ${\omega_{\infty}}^2=\omega_{\infty}-1$.
Noting that $\zeta_{\infty}\to1$ if $p\to\infty$,
$(\zeta_{\infty},\omega_{\infty})$ can be regarded as a solution to
\eqref{eq:algebraic_fig8_original} for a large $p$ .
\par
Now we calculate $V(\zeta_{\infty},\omega_{\infty})$.
Since
$\dfrac{\partial\,\tilde{V}}{\partial\,z}(\zeta_{\infty},\omega_{\infty})
\to0$
and
$\dfrac{\partial\,\tilde{V}}{\partial\,w}(\zeta_{\infty},\omega_{\infty})
\to0$
if $p\to\infty$,
\begin{align*}
  V(\zeta_{\infty},\omega_{\infty})
  &=
  -\Li_2(\zeta_{\infty}\omega_{\infty})
  +\Li_2\left(\dfrac{\zeta_{\infty}}{\omega_{\infty}}\right)
  +\dfrac{p}{4}\left(\log\zeta_{\infty}\right)^2
  -\log{\zeta_{\infty}}\log{\omega_{\infty}}
  \\
  &\xrightarrow[N\to\infty]{}
  -\Li_2\left(\exp\left(-\dfrac{\pi\sqrt{-1}}{3}\right)\right)
  +\Li_2\left(\exp\left(\dfrac{\pi\sqrt{-1}}{3}\right)\right)
  \\
  &=
   \pi^2
   \left\{
      \overline{B}_2\left(\dfrac{1}{6}\right)
     -\overline{B}_2\left(-\dfrac{1}{6}\right)
   \right\}
   +\sqrt{-1}
   \left\{
     \Lob\left(\dfrac{\pi}{3}\right)
    -\Lob\left(-\dfrac{\pi}{3}\right)
   \right\}
  \\
  &=
  \sqrt{-1}\times2\Lob\left(\dfrac{\pi}{3}\right).
\end{align*}
Here I used the fact that
\begin{equation*}
  \Li_2(\exp(\theta\sqrt{-1}))=
  \pi^2\overline{B}_2\left(\dfrac{\theta}{2\pi}\right)
  +\sqrt{-1}\Lob(\theta)
\end{equation*}
with $\Lob$ the Lobachevskij function and $\overline{B}_2$ the second
modified Bernoulli polynomial \cite[Proposition B]{Kirillov:dilog}:
\begin{equation*}
  \overline{B}_2(x):=
  -\dfrac{1}{\pi^2}\sum_{n=1}^{\infty}\dfrac{\cos(2n\pi x)}{n^2}.
\end{equation*}
Note that $\Lob$ is an odd function and $\overline{B}_2$ is an even function.
\par
Thus we can write
\begin{equation*}
  V(\zeta_{\infty},\omega_{\infty})
  =
  \CS(4_1)+\sqrt{-1}\Vol(4_1).
\end{equation*}
This suggests that the series of optimistic limits
$\displaystyle
  \left\{
    \olim_{N\to\infty}\dfrac{2\pi\sqrt{-1}\log\tau_{N}(M_p)}{N}
  \right\}_{p=0,1,\dots,\infty}$
goes to
$\displaystyle\lim_{N\to\infty}\dfrac{2\pi\sqrt{-1}\log{J_{N}(4_1)}}{N}$,
agreeing with the facts that
$\displaystyle\lim_{p\to\infty}\Vol(M_p)=\Vol(4_1)$
and that
$\displaystyle\lim_{p\to\infty}\CS(M_p)=\CS(4_1)$.
Note that $\CS(4_1)=0$ since $4_1$ is amphicheiral.
\begin{rem}
It was suggested by A.~Kricker to observe the $p\to\infty$ limit.
\end{rem}
\section{Volume Conjecture for closed three-manifolds}
Now I propose a very ambiguous conjecture.
\begin{conj}[Volume conjecture for closed three-manifolds]
  For any closed three-manifold $M$
  \begin{equation*}
    \olim_{N\to\infty}\dfrac{2\pi\sqrt{-1}\log\tau_N(M)}{N}
    =\CS(M)+\sqrt{-1}\Vol(M),
  \end{equation*}
where $\Vol(M):=v_3\Vert{M}\Vert$ with $\Vert{M}\Vert$ the simplicial volume
of $M$.
\end{conj}
A weaker but more precise conjecture is
\begin{conj}\label{conj:fig8}
For an integer $p$ put
\begin{equation*}
  V_p(z,w):=
  -\Li_2(zw)+\Li_2\left(\dfrac{z}{w}\right)+
  \dfrac{p}{4}(\log z)^2-\log z\log w.
\end{equation*}
Then there exists $(\zeta_p,\omega_p)$ such that
\begin{enumerate}
  \item
    $\dfrac{\partial\,V_p}{\partial\,z}(\zeta_p,\omega_p)
    =\dfrac{\partial\,V_p}{\partial\,w}(\zeta_p,\omega_p)=0$, and
  \item
    $V_p(\zeta_p,\omega_p)=\CS(M_p)+\sqrt{-1}\Vol(M_p)$.
\end{enumerate}
Here $M_p$ is the closed three-manifold obtained from the three-sphere by
$p$-surgery along the figure-eight knot.
\end{conj}
\begin{rem}
Conjecture~\ref{conj:fig8} is {\em numerically} true (up to 8 digits or so)
for $p=0,\pm1,\pm2,\pm3,\pm5,\pm6,\dots,\pm100$.
\end{rem}
\begin{rem}
Note that $V_p(1,w)=-\Li_2(w)+\Li_2\left(\dfrac{1}{w}\right)$ and this appears
in calculations about the figure-eight knot \cite[3.15]{Kashaev:LETMP97}.
More precisely $\omega:=\exp\left(-\dfrac{\pi\sqrt{-1}}{3}\right)$
satisfies
\begin{enumerate}
  \item $\dfrac{\partial\,V_p}{\partial\,w}(1,\omega)=0$, and
  \item $V_p(1,\omega)=\CS(4_1)+\sqrt{-1}\Vol(4_1)$.
\end{enumerate}
\end{rem}
Finally I raise some natural problems.
\begin{prob}\label{prob:fig8}
Can one generalize Conjecture~\ref{conj:fig8} to rational surgery?
Compare with \cite{Yoshida:INVEM85}.
\end{prob}
\begin{prob}
For any knot $K$ (or more generally link) and any {\em rational} number $p$,
does there exist a function $V_{p}$ as above?
\end{prob}
\begin{rem}
  Y.~Yokota told me that Conjecture~\ref{conj:fig8} and Problem~\ref{prob:fig8}
  can be solved by considering deformations of tetrahedron decomposition
  described in \cite{Yokota:volume}.
\end{rem}
\bibliography{mrabbrev,hitoshi}
\bibliographystyle{amsplain}
\end{document}